\documentclass[12pt]{amsart}

\DeclareMathOperator{\Br}{Br} \DeclareMathOperator{\Cor}{Cor}
\DeclareMathOperator{\h}{ht} \DeclareMathOperator{\Hom}{Hom}
\DeclareMathOperator{\Ind}{Ind} 

\usepackage{amscd}
\usepackage{amsfonts}
\usepackage{amssymb}
\usepackage{enumerate}

\begin{document}

\newtheorem*{theorem*}{Main Theorem}
\newtheorem*{corollary*}{Corollary}
\newtheorem{proposition}{Proposition}
\newtheorem{lemma}{Lemma}

\theoremstyle{definition}
\newtheorem*{remark*}{Remark}

\newcommand{\B}{\mathrm{B}}
\newcommand{\F}{\mathbb{F}}
\newcommand{\N}{\mathbb{N}}
\newcommand{\p}{\mathfrak{p}}
\newcommand{\q}{\mathfrak{q}}
\newcommand{\Q}{\mathbb{Q}}
\newcommand{\X}{\mathrm{X}}
\newcommand{\Xn}{\mathrm{X}_{\neq u}}
\newcommand{\Y}{\mathrm{Y}}
\newcommand{\Z}{\mathbb{Z}}

\parskip=10pt plus 2pt minus 2pt

\title[Brauer Groups of Genus Zero Extensions]
{Brauer Groups of Genus Zero Extensions of Number Fields}

\author[Jack Sonn]{Jack Sonn$^*$}
\address{Department of Mathematics, Technion---Israel
Institute of Technology, Haifa \ 32000 \ ISRAEL}
\thanks{$^*$Research supported by the Fund for
Promotion of Research at the Technion.}
\email{sonn@math.technion.ac.il}

\author[John Swallow]{John Swallow$^{**}$}
\address{Department of Mathematics, Davidson College,
Box 7046, Davidson, North Carolina \ 28035-7046 \ USA}
\thanks{$^{**}$Research supported in part by an International Research
Fellowship, awarded by the National Science Foundation
(INT--980199) and held at the Technion---Israel Institute of
Technology during 1998--1999, and a Young Investigator Grant from
the National Security Agency (MDA904-02-1-0061).}
\email{joswallow@davidson.edu}

\begin{abstract}
We determine the isomorphism class of the Brauer groups of certain
nonrational genus zero extensions of number fields.  In
particular, for all genus zero extensions $E$ of the rational
numbers $\Q$ that are split by $\Q(\sqrt{2})$, $\Br(E)\cong
\Br(\Q(t))$.
\end{abstract}

\date{January 29, 2003}

\maketitle

\section*{Introduction}

Let $K$ be a countable field, $\Br(K)$ its Brauer group. Then
$\Br(K)$ is a countable abelian torsion group; hence, as an
abstract group, it is completely determined by its Ulm invariants
$U_p(\lambda, \Br(K))$, where $p$ is a prime number and $\lambda$
is an ordinal number (see definition below).  With this
observation, Fein and Schacher initiated the investigation of the
Ulm invariants of algebraic function fields over global fields,
which culminated in the determination of all the Ulm invariants of
the Brauer group of a \emph{rational} function field $K$ in
finitely many variables over a global field $k$ (see \cite{FSS1}).
For nonrational function fields in one variable over a global
field $k$, the closest thing to a rational function field would be
a nonrational function field $K$ of genus zero over $k$. In
\cite{FSS2}, all the Ulm invariants $U_p(\lambda, \Br(K))$ for
such a field $K$ are determined, except for one (!), namely
$U_2(\omega2, \Br(K))$.  It turns out that all of the Ulm
invariants of $\Br(K)$, except for the missing one, coincide with
those of $\Br(k(t))$. The problem of this missing Ulm invariant
has remained open; in fact, until now, there was not even a single
example known of a nonrational genus zero function field over a
global field for which the missing Ulm invariant $U_2(\omega2,
\Br(K))$ had been computed.  In view of the fact that
$U_2(\omega2, \Br(k(t)))=0$, $\Br(K)\cong \Br(k(t))$ if and only
if $U_2(\omega2, \Br(K))=0$.

In this paper we determine the first known isomorphism classes of
Brauer groups of nonrational genus zero extensions $E/k$ of number
fields $k$. We prove the following
\begin{theorem*}
Let $k$ be a totally real number field and $l/k$ the quadratic
subfield of the cyclotomic $\Z_{2}$-extension $k^{cyc}/k$.

Suppose that $E$ is a genus zero extension of $k$ which is split
by $l$, and suppose that the Leopoldt Conjecture holds for the
field $l$ and the prime $2$.  Then $\Br(E)\cong \Br(k(t))$.
\end{theorem*}
Since the Leopoldt Conjecture (see \cite[\S 10.3]{NSW}) is known
to be true for all abelian number fields \cite{Br} (see also
\cite[Thms.~10.3.14 and 10.3.16]{NSW}), we have
\begin{corollary*}
    Suppose $k$ is a totally real abelian number field not
    containing $\sqrt{2}$.  Then if $E/k$ is a genus zero
    extension split by $k(\sqrt{2})$, then $\Br(E)\cong
    \Br(k(t))$.

    In particular, $\Br(E)\cong \Br(\Q(t))$ for all genus zero
    extensions $E$ of $\Q$ split by $\Q(\sqrt{2})$.
\end{corollary*}

We approach $\Br(E)$ via $\Br(El)$, where $l/k$ is a quadratic
extension such that $El$ is a rational function field. The
Auslander-Brumer-Faddeev theorem establishes an isomorphism
between $\Br(El)$ and a direct sum of $\Br(l)$ and character
groups of extensions of $l$, and we study an order two action on
$\Br(El)$ in terms of standard cohomological maps on these
summands.  Section~\ref{se:prelim} introduces this approach, and
section~\ref{se:actions} establishes properties of the order two
action.  Then in section~\ref{se:subgroup} we present a technical
analysis of heights of elements in the fixed subgroup of this
order two action.  Finally, in section~\ref{se:proof} we prove the
Main Theorem.

\section{Preliminaries}\label{se:prelim}

Let $k$ be a perfect field.  For any Galois field extension $K/k$,
let $G_{K/k}$ denote the Galois group.  Let $\bar K$ denote the
algebraic closure of $K$ and $G_{K}$ the absolute Galois group
$G_{\bar K/K}$.

We denote by $k^*$ the multiplicative group of $k$, by $\X(k)$ its
character group $H^1(G_k, \Q/\Z)$, and by $\Br(k)$ its Brauer
group $H^{2}(G_{k}, \bar k^{*})$.  Similarly, $\Br(K/k)$  denotes
the relative Brauer group of $K/k$, identified with
$H^{2}(G_{K/k}, K^{*})$.  All cohomology groups will be written
additively, and all modules will be left modules.

For an additive abelian torsion group $A$, we write $A_2$ for the
$2$-primary component of $A$.  The Ulm subgroups (for the prime 2)
of $A_2$ are defined for any ordinary $\lambda$ by $A_2(0)=A_2$,
$A_2(\lambda+1)=2A_2(\lambda)$, and for $\lambda$ a limit ordinal,
$A_2(\lambda)=\cap_{\lambda'<\lambda} A_2(\lambda')$.  The least
$\lambda$ such that $A_2(\lambda)=A_2(\lambda+1)$ is the Ulm
length of $A_2$.  Now let $P(\lambda)=\{\alpha\in A_2(\lambda) :
2\alpha=0\}$.  The Ulm invariant of $A_2$ at $\lambda$, denoted
$U_2(\lambda, A_2)$, is $[P(\lambda)/P(\lambda+1):\Z/2\Z]$.

We write $\h_{A_2}(\alpha)$ or $\h_A(\alpha)$ for the height of
$\alpha\in A_2$ in $A_2$, defined to be $\lambda$ such that
$\alpha\in A_2(\lambda)\setminus A_2(\lambda+1)$ if such a
$\lambda$ exists; otherwise, $\alpha\in A_2$ is divisible and we
write $\h_{A}(\alpha)=\infty$. We denote by $D(A)$ the divisible
subgroup of $A$.  For $\alpha\neq 0$, we will say that a
\emph{divisible tower over} $\alpha$ is a set
$\{\alpha_i\}_{i=0}^\infty$ satisfying $2^j\alpha_i= \alpha_{i-j}$
for $j\le i$ and $\alpha_0=\alpha$. Given nonzero $\alpha\in A_2$,
$\alpha\in D(A_2)$ if and only if there exists a divisible tower
over $\alpha$.

\subsection{The Genus Zero Extension $E$}\

A genus zero extension $E$ over a field $k$ is the quotient field
of $k[x,y]/\langle 1-cx^{2}-dy^{2} \rangle$ for $c, d\in k^*$. We
set $l=k(\sqrt{d})$.  The quaternion algebra $(c,d)$ is split in
$\Br(k)$ iff $E$ is isomorphic to a rational function field in one
variable over $k$. We will assume throughout that $(c,d)$ is not
split. Hence $l\not\subset E$.

Let $El$ denote the compositum of $l$ and $E$.  We determine an
element $u$ such that $El=l(u)$ as follows. Let $u = (1+\sqrt{d}y)
/ x$. Then $u(1-\sqrt{d}y)/x=c$, so that $c/u = (1-\sqrt{d}y)/ x$.
Moreover, $x=2/(u+c/u)$ and $y = (2u/(u+c/u)-1)/ \sqrt{d}$,
establishing that $El=l(u)$.

We denote by $\sigma$ the unique nontrivial element of $G_{l/k}$
and by $s$ the nontrivial element of $G_{El/E}$.  We may extend
$\sigma$ to an element of $G_k$, which by abuse of language we
also denote by $\sigma$, and we similarly denote by $s$ the
corresponding extension of $s$ to $G_{E\bar k/E}$.

In what follows $p$ will always be restricted to monic,
irreducible $p\in l[u]\setminus \{u\}$, and for each $p$ we fix a
root $a_p$ of $p$ for the remainder of the paper.  If
$\{a_{i}\}\subset \bar k^{*}$ are the roots of $p$, we set $\tilde
p\in l[X]$ to be the monic, irreducible polynomial with roots
$\{c/ \sigma(a_{i})\}$. Since $\sigma^2\in G_l$, $p\mapsto \tilde
p$ is an order two action.

\begin{lemma}\label{le:sact}
    The map $s$ acts triangularly on
    \begin{equation*}
        {\bar k}(u)^* = \bar k^* \times \langle u\rangle
        \times \coprod_{a\in \bar k^*} \langle u+a\rangle
    \end{equation*}
    by
    \begin{enumerate}
    \item $s(a)=\sigma(a)$, $a\in \bar k^*$;
    \item $s(u)=c\cdot \frac{1}{u}$; and
    \item $s(u+a) = \sigma(a) \cdot \frac{1}{u}
        \cdot (u+\frac{c}{\sigma(a)})$, $a\in \bar k^{*}$.
    \end{enumerate}
    Moreover, for all $p$, $s(p) =
    \sigma(p(0))\cdot u^{-\deg p}\cdot \tilde p(u)$.
\end{lemma}
For later use we define componentwise homomorphisms $s_{uu} \colon
\langle u\rangle \to \langle u\rangle$, $s_{u+a,u} \colon \langle
u+a\rangle\to \langle u\rangle$, and $s_{p\tilde p}\colon
\prod_{p(a)=0} \langle u+a\rangle\to \prod_{\tilde p(a)=0} \langle
u+a\rangle$ by $s_{uu}(u^e) = s_{u+a,u}((u+a)^e)=1/u^e$ and
\begin{equation*}
    s_{p\tilde p} \left( \prod_{p(a)=0} (u+a)^{e_a} \right) =
    \prod_{p(a)=0} \left(u+\frac{c}{\sigma(a)}\right)^{e_a} =
    \prod_{\tilde p(a)=0} (u+a)^{e_{c/\sigma^{-1}(a)}}.
\end{equation*}

\subsection{The Brauer Group of $E$}\label{se:brgrp}\

We have identified $\Br(E)$ with $H^{2}(G_{E}, \bar E^{*})$.  The
algebraic closure $\bar E$ of $E$ is identical to the algebraic
closure of $\bar k(u)$, and since by Tsen's theorem the Brauer
group of $\bar k(u)$ is trivial, we have that $H^2(G_{E}, \bar
E^{*})\cong H^{2}(G_{\bar k(u)/E}, \bar k(u)^{*})$.  Now every
element $\gamma\in G_{k}$ lifts to $\gamma' \in G_{\bar k(u)/E}$
by extending the automorphism trivially on $x$ and $y$. Inversely,
since $\bar k(u) = E\otimes_{k} \bar k$ and $k$ is algebraically
closed in $E$, any element $\gamma' \in G_{\bar k(u)/E}$ sends
$\bar k$ to $\bar k$. Therefore we may and do identify $G_{\bar
k(u)/E}$ with $G_{k}$, and we have $\Br(E)\cong H^{2}(G_{k}, \bar
k(u)^{*})$.

Now by Hilbert's Theorem 90, $H^1(G_l,\bar k(u)^*)$ is trivial.
Moreover, since $G_{l/k}$ is finite cyclic, $H^3(G_{l/k},l(u)^*)
\cong H^1(G_{l/k},l(u)^*)$, which is also trivial by Theorem 90.
The standard inflation-restriction five-term exact sequence
\cite[Prop.~1.6.6]{NSW} beginning with $H^2$ then begins
\begin{equation}\label{eq:brexactseq}
    0\to H^2(G_{l/k},l(u)^*)\to H^2(G_k,\bar k(u)^*)
    \xrightarrow{\phi} H^2(G_l,\bar k(u)^*)^{G_{l/k}}\to 0,
\end{equation}
or, equivalently,
\begin{equation*}
    0\to \Br(El/E) \to \Br(E)\to \Br(l(u))^{G_{l/k}} \to 0.
\end{equation*}

Now $\sigma \in G_k$ acts naturally on $\Br(l(u))=H^2(G_l,\bar
k(u)^*)$: given a 2-cocycle $h$, this action is $h^s(g_1, g_2) =
s(h(\sigma^{-1}g_1 \sigma, \sigma^{-1}g_2 \sigma))$, where $s$
acts on $\bar k(u)^*$ as above. We denote this action on
$\Br(l(u))$ by $s^*$, and the fixed group in \eqref{eq:brexactseq}
on the right is therefore $\Br(l(u))^{\langle s^*\rangle}$.

\subsection{The Auslander-Brumer-Faddeev Decomposition}\

Recall that the Auslander-Brumer-Faddeev theorem
(\cite[Prop.~4.1]{AB}, \cite[Thms. 15.2, 15.3]{F}) establishes an
isomorphism
\begin{equation}\label{eq:fund}
    \Br(l(u)) \cong \Br(l) \oplus \X(l) \oplus \left(\oplus_p
    \X(l(a_p))\right).
\end{equation}
We will need the particular isomorphisms contained in the proof of
this result, which we review as follows.

Let $G=G_l$ and for an arbitrary fixed $p$, $H=G_{l(a_p)}$. Let
$A= \Ind_{G}^{H}(\Z) = \Hom_H(G, \Z)$ be the $G$-module of
$H$-module homomorphisms from $G$ to the trivial $G$-module $\Z$;
$g\in G$ acts on $x\in \Ind_{G}^{H}(\Z)$ via $(gx)(g_1)= x(g_1
g)$. $A$ may be considered the set of functions from $G$ to $\Z$
defined on right cosets of $H$ in $G$. Define $\bar x(g_1)=
x(g_1^{-1})$ for $x\in A$.  Then $\bar x$ is defined on left
cosets of $H$ in $G$ and $\overline{gx}(g_1) = \bar x(g^{-1}g_1)$.

Consider $B = B(p) := \prod_{p(a)=0} \langle u-a\rangle \subset
\bar k(u)^*$. We claim that $A$ and $B$ are isomorphic $G$-modules
under the map $\iota\colon A\to B$ given by
\begin{equation*}
    \iota(x) = \prod_{\tau H\in G/H} (u-\tau(a_p))^{\bar x(\tau)}.
\end{equation*}
To check that $\iota$ respects $G$-action, we calculate
\begin{align*}
    \iota({gx}) &= \prod_{\tau H\in G/H}
    (u-\tau(a_p))^{\overline{gx}(\tau)} = \prod_{\tau H\in G/H}
    (u-\tau(a_p))^{\bar x(g^{-1}\tau)}
    \\
    &= \prod_{g\tau H\in G/H} (u-g\tau(a_p))^{\bar x(\tau)} = g
    (\prod_{\tau H\in G/H} (u-\tau(a_p))^{\bar x(\tau)}) =
    g\iota(x).
\end{align*}

The proof of the Auslander-Brumer-Faddeev theorem proceeds by
splitting the $G_l$-module $\bar k(u)^*$ of $H^2(G_l,\bar k(u)^*)$
as
\begin{equation}\label{eq:sp}
    \bar k(u)^* = \bar k^* \times \langle u\rangle \times
    \coprod_{p} \prod_{p(a)=0} \langle u-a\rangle,
\end{equation}
and then realizing the summands $\X(l)$ and $\X(l(a_p))$ with the
isomorphisms
\begin{equation}\label{eq:chiofl}
    \X(l) = H^1(G_l,\Q/\Z)\xrightarrow{\delta} H^2(G_l,\Z)
    \xrightarrow{\pi^*} H^2(G_l,\langle u\rangle)
\end{equation}
and
\begin{equation}\label{eq:chioflap}
    \X(l(a_p)) = H^1(H, \Q/\Z) \xrightarrow{\delta} H^2(H, \Z)
    \xrightarrow{sh} H^2(G,A) \xrightarrow{\iota^*} H^2(G, B).
\end{equation}
Here $\delta$ denotes the standard coboundary map, $sh$ the map of
Shapiro's Lemma, and $\pi\colon \Z\to \langle u\rangle$ the
homomorphism $\pi(e)=u^e$ of trivial $G_l$-modules $\Z$ and
$\langle u\rangle$.

We identify $\Br(l(u))$ with the decomposition \eqref{eq:fund} and
denote an arbitrary element of this group by $\beta+\chi_u+ \sum
\chi_p$ or $\beta \oplus\chi_u \oplus \sum \chi_p$.

By Lemma~\ref{le:sact}, since $s$ acts triangularly on $\bar
k(u)^*$ on the factors in \eqref{eq:sp} and $s_{p\tilde p}$ sends
$B(p)$ to $B(\tilde p)$, we have

\begin{lemma}\label{le:sstaract}
The map $s^*$ acts triangularly on $\Br(l(u))$. In particular, for
arbitrary $\beta \oplus \chi_u \oplus \sum\chi_p \in \Br(l(u))$,
\begin{equation}\label{eq:sstar}
\begin{split}
    s^*(\beta \oplus \chi_u \oplus \sum \chi_p) &= \left(
    s_{11}^*(\beta) + s_{u1}^*(\chi_u) + \sum s_{p1}^*(\chi_p)
    \right) \oplus \\ &\phantom{=} \ \left( s_{uu}^*(\chi_u) +
    \sum s_{pu}^*(\chi_p) \right) \oplus \\ &\phantom{=} \ \sum
    s_{p\tilde p}(\chi_p),
\end{split}
\end{equation}
where we denote the component parts of $s^*$ as follows:
\begin{enumerate}[\indent \indent \indent (a)\indent]
\item $s_{p\tilde p}^*\colon \X(l(a_p))\to \X(l(a_{\tilde p}))$;
\item $s_{pu}^* \colon \X(l(a_p))\to \X(l)$; \item $s_{uu}^*\colon
\X(l)\to \X(l)$; \item $s_{p1}^* \colon \X(l(a_p))\to \Br(l)$;
\item $s_{u1}^* \colon \X(l)\to \Br(l)$; and \item $s_{11}^*
\colon \Br(l)\to \Br(l)$.
\end{enumerate}
\end{lemma}

For later use we further define $s_u^*=\sum_p s_{pu}^*$ and
$s_1^*=s_{u1}^*+\sum_p s_{p1}^*$.

\section{Decomposing $s^*$ on $\Br(l(u))$}\label{se:actions}

In this section we study the component functions of $s^*$,
comparing them to natural Galois actions on $\Br(l)$ and $\X(l)$.
We denote by $\sigma$ the natural Galois action of $\sigma$ on
$\Br(l)\cong H^2(G_l, \bar k^*)$ and on $\X(l)\cong H^1(G_l,
\Q/\Z)$. Note that in $H^1(G_l, \Q/\Z)$, $\Q/\Z$ is a trivial
$G_k$-module.

\subsection{$s_{11}^*$ on $\Br(l)$ and $\Br(l_\p)$}\

Since $s$ acts on $\bar k^*$ via $\sigma$, the action of
$s_{11}^*$ on $\Br(l)$ is the natural Galois action:
\begin{equation*}
    s_{11}^*(\beta) = \sigma(\beta), \qquad \beta\in \Br(l).
\end{equation*}
We now determine the action of $\sigma$ on the local invariants
$b_\p$ of an element $b$ of $\Br(l)$ with the following general
result.

\begin{proposition}\label{pr:action}
    Let $K/k$ be a finite Galois extension of global fields.
    Suppose $g\in G_{K/k}$ and $\p$ be a prime of $K$, and let
    $b_\p$ denote the local invariant of $b\in \Br(K)$ at $\p$.
    Then
    \begin{equation}\label{eq:action}
        (g(b))_\p = b_{g^{-1}(\p)}.
    \end{equation}
\end{proposition}

\begin{proof}
    For an arbitrary field $K$, a $K$-algebra $A$ is an
    associative ring $A$ with unity together with an embedding
    $\alpha\colon K\hookrightarrow A$ of $K$ into the center of
    $A$.  Thus a $K$-algebra should be considered a pair
    $(A,\alpha)$.  Now let $(B,\beta)$ be a second $K$-algebra.
    Then the tensor product
    \begin{equation*}
        (C, \gamma)=(A, \alpha)\otimes_K(B, \beta)
    \end{equation*}
    is generated by elements $a\otimes b$, $a\in A$, $b \in B$,
    satisfying $\alpha(r)a\otimes b=a\otimes \beta(r)b$, $r \in
    K$, and where $\gamma=\alpha \odot \beta \colon K
    \hookrightarrow C$ is defined by $\gamma(r) = \alpha(r)
    \otimes 1_B = 1_A \otimes \beta(r)$.

    Let $g$ be an automorphism of $K$. An action of $g$ on $(A,
    \alpha)$ may be defined by $(A, \alpha)^{g} := (A,
    \alpha^{g})$, where $\alpha^{g}(r)=\alpha g(r)$ for $r\in K$.
    Note that as rings, $(A,\alpha)$ and $(A,\alpha)^{g}$ are
    isomorphic. It follows that $\gamma^{g}(r)=\gamma g(r)=\alpha
    g(r)\otimes 1 =1\otimes\beta g(r)$, \emph{i.e.},
    $\alpha^{g}(r)\otimes 1 = 1\otimes \beta^{g}(r)$, so
    \begin{equation}\label{eq:tensoract}
        \gamma^{g}= (\alpha\odot\beta)^{g} = \alpha^{g}
        \odot \beta^{g}.
    \end{equation}

    We now specialize to the case where $(A, \alpha)$ is a
    finite-dimensional central simple $K$-algebra and $B$ is a
    complete discrete valued field $K_\p$, where $\p$ will denote
    the valuation of $K_\p$ as well as the valuation induced on
    $K$ by the embedding $\beta$ of $K$ into $K_\p$. By
    \eqref{eq:tensoract},
    \begin{equation*}
        (\alpha^{{g}^{-1}}\odot \beta)^{g} =
        \alpha\odot\beta^{g}.
    \end{equation*}
    Note again that, as rings,  $(C, (\alpha^{{g}^{-1}} \odot
    \beta)^{g})$ and $(C, \alpha^{{g}^{-1}} \odot \beta)$ are
    isomorphic. Moreover, $(C, (\alpha^{{g}^{-1}} \odot
    \beta)^{g})$ and $(C, \alpha\odot\beta^{g})$ both represent
    finite-dimensional central simple algebras over complete
    discrete valued fields.

    Now assume the hypotheses of the Proposition.  By a lemma of
    Janusz \cite[Lemma, p.~385]{Ja}, $(C, (\alpha^{{g}^{-1}} \odot
    \beta)^{g})$ and $(C, \alpha^{{g}^{-1}}\odot\beta)$ have the
    same local invariant. It follows that $(C, \alpha^{{g}^{-1}}
    \odot \beta)$ and $(C, \alpha \odot \beta^{g})$ have the same
    local invariant.
\end{proof}

\begin{remark*}
    In the proof of Janusz' lemma, the assertion is made that a
    field isomorphism of $p$-adic local fields is an isomorphism
    of $p$-adic local fields, \emph{i.e.}, that prime elements are
    mapped to prime elements.  To prove this statement, the
    essential idea is the same as in the proof that the only
    algebraic automorphism of $\Q_p$ is the identity. Suppose
    $L_i$, $i=1, 2$ are finite extensions of $\Q_p$ with
    ramification indices $e_i$ and normalized valuations $v_i$,
    $f\colon L_1\to L_2$ is a field isomorphism, and $\pi$ is a
    prime element of $L_1$.  Then $\pi^{e_1}=pu$ with $u\in L_1$ a
    unit.  Write $f(\pi)^{e_1} = f(pu) = f(p) f(u) = p f(u)$.  Now
    if $f(u)$ is a unit in $L_2$, we are done: since
    $e_1v_2(f(\pi))=v_2(p)=e_2$, $e_1 \mid e_2$ and by symmetry
    $e_1=e_2$ and $v_2(f(\pi))=1$. That $f$ maps units to units
    follows from the fact that the units of a $p$-adic field may
    be characterized algebraically, as the only elements that have
    $n$th roots in the field for infinitely many $n$.
\end{remark*}

\subsection{$s_{uu}^*$ on $\X(l)$} \

\begin{proposition}\label{pr:suuact}
    In the decomposition \eqref{eq:fund},
    \begin{equation}\label{eq:suuact}
        s_{uu}^*(\chi) = -\sigma(\chi), \quad \chi\in \X(l),
    \end{equation}
    where $\sigma$ denotes the natural Galois action on
    $\X(l)=H^{1}(G_{l},\Q/\Z)$.
\end{proposition}

\begin{proof}
    The homomorphism $s_{uu}^*$ acts on a 2-cocycle $f\in
    H^{2}(G_{l},\langle u\rangle)$ by
    \begin{equation*}
        f(g_1, g_2) \mapsto s_{uu}(f(\sigma^{-1}g_1\sigma,
        \sigma^{-1}g_2\sigma)) = -f(\sigma^{-1}g_1\sigma,
        \sigma^{-1}g_2\sigma),
    \end{equation*}
    since the image of $f$ lies in $\langle u\rangle$ and $s_{uu}$
    is the inversion map.  Now $\X(l)$ appears in \eqref{eq:fund}
    under the isomorphisms in \eqref{eq:chiofl}.  Using the facts
    that in the natural Galois actions on $H^1(G_l, \Q/\Z)$ and
    $H^2(G_l, \Z)$, the action of $\sigma$ on the modules $\Q/\Z$
    and $\Z$ is trivial, and that $\sigma$ commutes with coboundary
    maps $\delta$ \cite[Prop.~1.5.2]{NSW},  a direct calculation
    yields
    \begin{equation*}
        \delta^{-1} \circ (\pi^*)^{-1} \circ s_{uu}^* \circ \pi^*
        \circ \delta = -\sigma.
    \end{equation*}
\end{proof}

\subsection{$s_{pu}^*$ on $\X(l(a_p))$} \

\begin{proposition}\label{pr:spuact}
    In the decomposition \eqref{eq:fund},
    \begin{equation}\label{eq:spuact}
        s_{pu}^*(\chi) = -\sigma(\Cor_{l(a)/l} (\chi_p)),
        \quad \chi_p\in \X(l(a_p)),
    \end{equation}
    where $\sigma$ denotes the natural Galois action on
    $\X(l)=H^{1}(G_{l},\Q/\Z)$.
\end{proposition}

\begin{proof}
    Let $B=\prod_{p(a)=0} \langle u-a\rangle$. The homomorphism
    $s_{pu}^*$ acts on a 2-cocycle $f\in Z^{2}(G_{l}, B)$ by
    \begin{equation*}
        f(g_1, g_2) \mapsto (\prod s_{u+a,u})
        (f(\sigma^{-1}g_1\sigma,
        \sigma^{-1}g_2\sigma)),
    \end{equation*}
    where for each $a$, $s_{u+a,u}(u+a)=\frac{1}{u}$. Further let
    $\rho\colon B\to \langle u\rangle$ be defined by
    \begin{equation*}
        \rho\left( \prod_{p(a)=0} (u-a)^{e_a}\right)=
        u^{(\sum_{p(a)=0} e_a)},
    \end{equation*}
    and let $\rho^*$ be the induced map $\rho^*\colon H^2(G_l,
    B)\to H^2(G_l, \langle u\rangle)$ of cohomology. Observe that
    $s_{pu}=s_{uu}\circ \rho$, hence $s_{pu}^* = s_{uu}^* \circ
    \rho^*$.

    Keeping the notation for $G$, $H$, and $A$ as in
    \eqref{eq:chioflap}, define a $G$-module homomorphism
    $\nu\colon A \to \Z$ by
    \begin{equation*}
        \nu(x)=\sum_{\tau H\in G/H} \tau x(\tau^{-1})
        =\sum_{\tau H\in G/H} x(\tau^{-1}).
    \end{equation*}
    We calculate
    \begin{align*}
        \rho(\iota(x)) &= \rho\left(\prod_{\tau H\in G/H}
        (u-\tau(a_p))^{\bar x(\tau)}\right) = \prod_{\tau H\in
        G/H} u^{\bar x(\tau)} \\ &= u^{(\sum_{\tau H\in G/H}
        x(\tau^{-1}))} = \pi\left(\sum_{\tau H\in G/H}
        x(\tau^{-1})\right) = \pi(\nu(x)).
    \end{align*}
    Hence we have a commutative diagram
    \begin{equation*}
    \begin{CD}
        H^2(G, A) @>{\iota^*}>> H^2(G, B)\\
        @V{\nu^*}VV @VV{\rho^*}V \\
        H^2(G, \Z) @>{\pi^*}>>
        H^2(G, \langle u\rangle)\\
    \end{CD}
    \end{equation*}
    Now $\nu^*\circ sh$ is the corestriction map from
    $H^2(G_{l(a_p)},\Z)\to H^2(G_l,\Z)$ \cite[Prop.~1.6.4]{NSW}.
    Moreover, the coboundary map $\delta$ from $H^1$ to $H^2$
    commutes with the corestriction map \cite[Prop.~1.5.2]{NSW}.
    Hence under the isomorphisms in \eqref{eq:chiofl} and
    \eqref{eq:chioflap},
    \begin{align*}
        \delta^{-1} \circ (\pi^*)^{-1} \circ
        \rho^* \circ \iota^* \circ sh \circ \delta & =
        \delta^{-1} \circ \nu^* \circ sh \circ \delta \\
        &= \delta^{-1} \circ (\Cor_{H^2}) \circ \delta =
        \Cor_{H^1}.
    \end{align*}
    Therefore in the decomposition \eqref{eq:fund}, $\rho^*$ is
    the corestriction map, and since $s_{pu}^*=s_{uu}^*\circ
    \rho^*$, by \eqref{eq:suuact}, $s_{pu}^*=-\sigma
    \Cor_{l(a)/l}$.
\end{proof}

\subsection{$s_{p\tilde p}^*$ on $\X(l(a_p))$}\

Let $\Cor\colon \oplus \X(l(a_p))\to \X(l)$ be the sum $\sum
\Cor_{l(a_p)/l}$ of the corestrictions on each summand.

\begin{proposition}\label{pr:sppact}
    In the decomposition \eqref{eq:fund},
    \begin{equation}\label{eq:corsigma}
        \Cor_{l(a_{\tilde p})/l}(s_{p\tilde p}^*(\chi_p)) =
        \sigma(\Cor_{l(a)/l} (\chi_p)), \quad \chi_p\in
        \X(l(a_p))
    \end{equation}
    and if $p=\tilde p$, then
    \begin{equation}\label{eq:spppistp}
        s_{pp}^*(\chi_p) = \tilde \sigma(\chi_p), \quad
        \chi_p\in \X(l(a_p)),
    \end{equation}
    where $\tilde\sigma=\sigma\tau$ for a $\tau\in G_l$ such
    that $\tilde\sigma(a_p)=c/a_p$.

    If for all $p$, $s_{p\tilde p}^*(\chi_p)=
    \chi_{\tilde p}$, then
    \begin{equation}\label{eq:corsum}
        \sigma(\Cor \sum \chi_p) = \Cor \sum \chi_p
    \end{equation}
    and if additionally $\chi_u\in \X(l)$ such that
    $(1+\sigma)\chi_u = - \Cor \sum \chi_p$, then
    \begin{equation}\label{eq:corsum2}
        \sigma(s_1^*(\chi_u+\sum \chi_p)) =
        -s_1^*(\chi_u+\sum \chi_p).
    \end{equation}
\end{proposition}

\begin{proof}
    Equations \eqref{eq:corsigma}, \eqref{eq:corsum}, and
    \eqref{eq:corsum2} each follow from the fact that $s^*$ has
    order 2. For the first result, let $\chi_p\in \X(l(a_p))$ be
    arbitrary and set $\chi_{\tilde p}=s_{p\tilde p}^*(\chi_p)$.
    Then by applying \eqref{eq:sstar} twice, we see that the
    component of $s^*(s^*(\chi_p))=\chi_p$ in $\X(l)$ is
    \begin{equation*}
        s_{uu}^*s_{pu}^*(\chi_p)+s_{\tilde pu}^*(\chi_{\tilde p}),
    \end{equation*}
    which, by \eqref{eq:suuact} and \eqref{eq:spuact}, is
    $\Cor \chi_p -\sigma \Cor \chi_{\tilde p}$.  Then since this
    component must be trivial, the first result follows.
    A similar argument establishes \eqref{eq:corsum} and
    \eqref{eq:corsum2}.

    For \eqref{eq:spppistp}, first observe that if $p=\tilde p$,
    then the action $a\mapsto c/\sigma(a)$ permutes the roots of
    $p$.  Hence $c/\sigma(a)=a_p$ for some root $a$.  Since $G_l$
    is transitive on the roots of $p$, $a=\tau(a_p)$ for some
    $\tau \in G_l$.  Hence $a_p=c/\sigma\tau(a_p)$.  Let $\tilde
    \sigma = \sigma \tau$.

    Keeping the notation for $G$, $H$, $A$, and $B$ as in
    \eqref{eq:chioflap}, observe that $s_{pp}^*$ acts on
    $H^2(G,B)$ by sending an arbitrary 2-cocycle $f(g_1, g_2)$
    to
    \begin{equation*}
        s_{pp}(f(\sigma^{-1}g_1 \sigma, \sigma^{-1}g_2
        \sigma)).
    \end{equation*}

    Now $\tau^*$ acts on $H^2(G, B)$ trivially
    \cite[Prop.~1.6.2]{NSW}, so $s_{pp}^*\circ \tau^*=s_{pp}^*$.
    However, direct calculation shows that $s_{pp}^*\circ \tau^*$
    sends $f(g_1, g_2)$ to $(s_{pp}\circ
    \tau)(f({\tilde\sigma}^{-1}g_1 \tilde\sigma, \
    {\tilde\sigma}^{-1}g_2 \tilde\sigma))$, where
    \begin{align*}
        (s_{pp}\circ \tau)\left(\prod_{p(a)=0} (u-a)^{e_a}\right)
        &= \prod \left(u-\frac{c}{\sigma \tau(a)}\right)^{e_a} \\
        &= (u-a_p)^{e_{a_p}} \prod_{a\neq a_p}
        \left(u-\frac{c}{\sigma\tau(a)}\right)^{e_a}.
    \end{align*}
    Hence $s_{pp}\circ \tau$ acts trivially on $\langle
    u-a_p\rangle \subset B$, and consequently $\iota^{-1} \circ
    (s_{pp}\circ \tau) \circ \iota$ acts trivially on the values
    $x(1)$, $x\in A$.

    Recall that $(sh)^{-1}$ is induced by the module map $A\to \Z$
    given by $x\mapsto x(1)=x(H)$ \cite[Prop.~1.6.3]{NSW}; hence
    $(sh)^{-1} \circ (\iota^*)^{-1} \circ s_{pp}^* \circ \tau^*
    \circ \iota^* \circ sh$ acts on $H^2(H, \Z)$ by sending a
    2-cocycle $f(g_1, g_2)$ to $f({\tilde \sigma}^{-1}g_1
    \tilde\sigma, \ {\tilde \sigma}^{-1}g_2 \tilde\sigma)$. As a
    result,
    \begin{equation*}
        {\tilde \sigma} = sh^{-1} \circ (\iota^*)^{-1} \circ
        s_{pp}^* \circ \iota^* \circ sh
    \end{equation*}
    and since $\tilde\sigma$ commutes with coboundary maps
    \cite[Prop.~1.5.2]{NSW}, $s_{pp}^*(\chi_p) =
    \tilde\sigma(\chi_p)$.
\end{proof}

\section{The Fixed Subgroup $\Br(l(u))^{\langle
s^*\rangle}$}\label{se:subgroup}

\subsection{Notation}\

Let $\B^G$ denote the fixed subgroup $\Br(l(u))_2^{\langle s^*
\rangle}$ of $\Br(l(u))_2$.  By \eqref{eq:sstar} and
\eqref{eq:corsum}, a direct calculation shows that this group
consists of precisely those elements satisfying

\begin{enumerate}[\indent (i)]
    \item \label{cond1} $s_{p\tilde p}^* \chi_p = \chi_{\tilde
    p}$, $\forall p$;
    \item \label{cond2} $(1+\sigma)\chi_u = - \Cor \sum
    \chi_p$; and
    \item \label{cond3} $(1-\sigma)\beta = s_1^*(\chi_u+\sum
    \chi_p)$.
\end{enumerate}

Let $\Xn=\oplus_p \X(l(a_p))_2$ and $\X=\X(l)_2 \oplus \Xn$. We
further define $\Xn^G$ and $\X^G$ as follows:
\begin{align*}
    \Xn^G &= \left\{ \sum \chi_p \in \Xn \colon s_{p\tilde
    p}^*(\chi_p) = \chi_{\tilde p} \ \forall p\right\}; \\
    \X^G &= \left\{ \chi = \chi_u + \sum \chi_p \in \X \colon
    (1+\sigma)\chi_u = - \Cor \sum \chi_p, \right. \\
    &\phantom{=} \left.
    \ \ s_{p\tilde p}^* \chi_p = \chi_{\tilde p} \ \forall p
    \right\}.
\end{align*}

Now assume that $k$, and so also $l$, is a number field.  Then
$D(\X(l)_2)$ is of finite rank \cite[Thm.~11.1.2]{NSW} and its
dual $\Hom(D(\X(l)_2), \Q/\Z)$ is a free $\Z_2$-module of finite
rank.  Applying \cite[Thm.~74.3]{CR} and passing to the dual
again, we may then decompose $D(\X(l)_2)$ as
\begin{equation}\label{eq:dchil2}
    D(\X(l)_2) = I\oplus N\oplus P,
\end{equation}
where each of $I$, $N$, and $P$ is a finite direct sum of
$\Q_2/\Z_2$ summands and the natural Galois action $\sigma$ is
trivial on $I$, negation on $N$, and permutation of pairs of
$\Q_2/\Z_2$ summands on $P$. Then $D(\X(l)_2)^{\langle
\sigma\rangle}=I\oplus {}_2N\oplus P^{\langle \sigma\rangle}$.

\begin{remark*}
A sketch of an elementary proof of \eqref{eq:dchil2} is as
follows. Let $M$ be a direct sum of finitely many copies of
$\Q_2/\Z_2$, and let $\sigma$ be an order 2 action on $M$. Denote
by $V$ the $\F_2$-vector space given by the 2-torsion ${}_2M$ of
$M$. By linear algebra $V$ decomposes into $V_I\oplus V_P \oplus
\sigma(V_P)$.  Let $W_1$ be a complement of $(V_P+\sigma(V_P))\cap
(1+\sigma)M$ in $V\cap (1+\sigma)M$ and $W_{-1}$ a complement of
$(V_P+\sigma(V_P))\cap (1-\sigma)M$ in $V\cap (1-\sigma)M$.  Using
the fact that $M=2M=(1+\sigma)M+(1-\sigma)M$, one shows that
$V=V_1 \oplus V_{-1} \oplus V_P \oplus \sigma(V_P)$.  Now since
$M$ is divisible, the homomorphic image $(1+\sigma)M$ is
divisible; we may therefore construct a divisible tower inside
$(1+\sigma)M$ over each element of a basis of $W_1$.  Denote by
$I$ the subgroup generated by these towers. Similarly denote by
$N$ the subgroup generated by towers inside $(1-\sigma)M$
constructed over each element of a basis of $W_{-1}$.  Finally
construct towers over each basis element of $V_P$ and let $P$ be
the subgroup generated by these towers and their
$\sigma$-conjugates.  Observe that a sum $\sum M_i$ of subgroups
$M_i$ of $M$ is direct if and only if $\sum (M_i\cap V)=\oplus
(M_i\cap V)$. Moreover, if each $M_i$ is divisible and $\sum
\dim_{\F_2} (M_i\cap V) = \dim_{\F_2} V$, then $M=\oplus M_i$.
Hence $M=I\oplus N\oplus P$.
\end{remark*}

Each $w\in D(\X(l)_2)^{\langle \sigma\rangle}$ of order 2
represents a cyclic extension of $l$ of degree 2 which is fixed
under the action of $\sigma$.  By Kummer theory $w=l(\sqrt{e})$
for $e\in l^*\setminus {l^{*2}}$, and since $w$ is fixed by
$\sigma$, $l(\sqrt{e})$ is Galois over $k$. The group
$G_{l(\sqrt{e})/k}$ is isomorphic either to $\Z/4\Z$ or to $\Z/2\Z
\oplus \Z/2\Z$. One may choose a representative in $k^*$ for the
class of $e$ in $l^*/{l^{*2}}$ if and only if $G_{l(\sqrt{e})/k}
\cong \Z/2\Z \oplus \Z/2\Z$. In this case we write
$w=l(\sqrt{e})$, $e\in k^*\setminus k^{*2}$. We also write that
$w=0$ satisfies $w=l(\sqrt{1})$.

Define $W$ to be the exponent 2 subgroup of $D(\X(l)_2)^{\langle
\sigma\rangle}$ consisting of extensions with Klein 4-group over
$k$, together with the identity element:
\begin{equation*}
    W=\{ w\in D(\X(l)_2)^{\langle \sigma\rangle}: \vert w\vert \le
    2, \ w=l(\sqrt{e}), \ e\in k^*\}.
\end{equation*}
For an element $w\in D(\X(l)_2)$, write $w_I$ for the component of
$w$ in $I$, and set $W_I = \{ w_I\in I : w\in W\}$.

Now define a function $\lambda$ from $W_I$ to the set of ordinals
as follows:
\begin{equation*}
    \lambda_w = \sup \left\{ \h_{\Xn^G}\left(\sum \chi_p\right) :
    \left\vert\sum \chi_p \right\vert = 2, \ \ \left(\Cor \sum
    \chi_p\right)_I = w\right\}.
\end{equation*}
We will repeatedly use the fact that the Ulm length of $\X(K)$ for
$K$ a number field is less than $\omega 2$ \cite[Thm.~1]{FS}, and
hence that the Ulm lengths of $\X_{\neq u}$ and $\X$ are at most
$\omega 2$. Here, this result implies that for any $w\in W_I$ we
have that $\lambda_w \le \omega 2$ or $\lambda_w=\infty$.

\subsection{Reducing to $\X^G$}\

\begin{proposition}\label{pr:equiv1}
    Suppose $k$ is a number field.  Then there exists an order 2
    element in $\B^G$ of height $\omega 2$ if and only if there
    exists an order 2 element in $\X^G$ of height $\omega 2$.
\end{proposition}

\begin{lemma}\label{le:bggei}
    Suppose $k$ is a number field and that $\alpha=\beta + \chi\in
    \B^G$ has $\beta\in D(\Br(l))$ and $\h_{\X^G}(\chi)\ge i+1$ for
    some $i\in \N_{\ge 0}$.  Then $\h_{\B^G}(\alpha)\ge i$.
\end{lemma}

\begin{proof}
    We are given that there exists $\chi_{i+1}\in \X^G$ with
    $2^{i+1}\chi_{i+1}=\chi$.  Let $\chi_i=2\chi_{i+1}$.  By
    specifying the element by invariants at each place $\p$ of
    $l$, we will construct $\beta_i\in \Br(l)$ such that
    $2^i\beta_i=\beta$ and $(1-\sigma)\beta_i = s_1^*\chi_i$.
    Then $\alpha_i = \beta_i+\chi_i$ satisfies $\alpha_i\in \B^G$
    and $2^i\alpha_i=\alpha$.

    Consider places $\p_1, \p_2$ permuted by $\sigma$.  Let
    $x=\beta_{\p_1}$, $y=\beta_{\p_2}$.  Since $\alpha\in \B^G$,
    $(1-\sigma)\beta = s_1^*\chi$.  Then, by \eqref{eq:action},
    $(s_1^*\chi)_{\p_1} = ((1-\sigma)\beta)_{\p_1} =x-y$.  Let
    $z=(s_1^* \chi_{i+1})_{\p_1}$. Since $s_1^*$ is a
    homomorphism, $2^{i+1}z=x-y$. Set $(\beta_i)_{\p_1} = 2z +
    \frac{1}{2^i}y$ and $(\beta_i)_{\p_2} = \frac{1}{2^i} y$. Then
    $(2^i\beta_i)_{\p_1} =2^{i+1}z+y=x-y+y=x$ and
    $(2^i\beta_i)_{\p_2}=y$.  Now by \eqref{eq:action}
    \begin{equation*}
        ((1-\sigma)\beta_i)_{\p_1} = 2z + \frac{1}{2^i}y -
       \frac{1}{2^i}y = 2z = (s_1^*\chi_i)_{\p_1}.
    \end{equation*}
    Moreover, using \eqref{eq:action} and, by \eqref{eq:corsum2},
    that $\sigma(s_1^*(\chi)) = -s_1^*(\chi)$ for $\chi\in \X^G$,
    we obtain
    \begin{equation*}
        ((1-\sigma)\beta_i)_{\p_2} =
        -2z = (s_1^*\chi_i)_{\p_2}.
    \end{equation*}
    For all pairs $\p_j$, $j=1,2$ considered, then, we have
    $(2^i\beta_i)_{\p_j}=\beta_{\p_j}$ and
    $((1-\sigma)\beta_i)_{\p_j} = (s_1^*\chi_i)_{\p_j}$.

    Now consider archimedean, inert, or ramified places $\p$. By
    \eqref{eq:action}, for any $\gamma\in \Br(l)$,
    $((1-\sigma)\gamma)_{\p}=0$ at each such $\p$. Moreover, by
    \eqref{eq:corsum2}, for any $\delta\in \X^G$,
    $\sigma(s_1^*(\delta)) = -s_1^*(\delta)$. Hence
    $(s_1^*(\delta))_\p \in \{0,1/2\}$. If $\delta$ is divisible
    by $2$ in $\X^G$, then since $s_1^*$ is a homomorphism we have
    that $(s_1^*(\delta))_\p=0$ at any such $\p$. Hence the
    condition $((1-\sigma)\beta_i)_{\p}=(s_1^*\chi_i)_{\p}$ is
    always satisfied for any choice of $(\beta_i)_{\p}$.  We
    choose $(\beta_i)_{\p}$ for these $\p$ as follows.

    Since $\beta\in D(\Br(l))$, $\beta_\p=0$ at any archimedean
    $\p$. We define $(\beta_i)_\p = 0$ at all archimedean $\p$.
    Now consider the other $\p$, which are inert or ramified over
    $k$. Let $\q$ denote some such place with $\beta_{\q}=0$. Now
    for all such $\p\neq \q$, define $(\beta_i)_\p =
    \frac{1}{2^i}\beta_\p$ at such $\p$.  Then for all these $\p$
    considered, we have $(2^i\beta_i)_\p=\beta_\p$ and
    $((1-\sigma)\beta_i)_\p = (s_1^*\chi_i)_\p$. Now $\sum_{\p}
    \beta_{\p} = \sum_{\p\neq \q} \beta_{\p} = 0$ since
    $\beta_{\q}=0$.  Moreover, $2^i(\beta_i)_{\p} = \beta_{\p}$
    for all $\p\neq \q$.  Hence $2^i(\sum_{\p \neq \q}
    (\beta_i)_{\p}) = 0$.  Set $(\beta_i)_{\q} = - \sum_{\p \neq
    \q} (\beta_i)_{\p}$. Then $2^i(\beta_i)_{\q} = 0 = \beta_{\q}$
    and $\sum_{\p} (\beta_i)_{\p} = 0$.  Hence $\{\beta_{\p}\}$
    defines an element of $\Br(l)_2$ satisfying $2^i\beta_i=\beta$
    and $(1-\sigma)\beta_i = s_1^*\chi_i$ and we are done.
\end{proof}

\begin{lemma}\label{le:doublele}
    Suppose $k$ is a number field, and let $\alpha=\beta+\chi \in
    \B^G$. Suppose that $\h_{\B^G}(\alpha)>0$.  Then
    $\h_{\X^G}(\chi)-1\le \h_{\B^G}(\alpha) \le \h_{\X^G}(\chi)$.
\end{lemma}

\begin{proof}
    By projection from $\B^G$ to $\X^G$ we see that
    $\h_{\B^G}(\alpha)\le \h_{\X^G}(\chi)$.

    Suppose that $h=\h_{\X^G}(\chi)$ lies in $\N$.  Since
    $\h_{\B^G} (\alpha)>0$, $\beta_{\p}=0$ at all archimedean $\p$
    and $\beta\in D(\Br(l))$. By Lemma~\ref{le:bggei},
    $\h_{\B^G}(\alpha) \ge \h_{\X^G}(\chi)-1.$ Now suppose that
    $h=\omega$.  In this case Lemma~\ref{le:bggei} shows that in
    fact $\h_{\B^G}(\alpha) = \h_{\X^G}(\chi).$  Since the Ulm
    length of $\X^G$ is less than or equal to $\omega 2$
    \cite[Thm.~1]{FS}, we are left with the cases $h=\omega+n$ for
    some $n$, $h= \omega 2$, and $h=\infty$.

    For each $i\ge 1$ such that $h\ge \omega + i$ we do the
    following. Let $\chi_i\in \X^G$ be such that $2^i\chi_i=\chi$
    and for each $j \in \N$ let $\delta_j\in \X^G$ satisfy
    $2^j\delta_j = \chi_i$. The proof of Lemma~\ref{le:bggei}
    shows that we may find $\beta_{i-1}$ such that
    $\alpha_{i-1}=\beta_{i-1}+\chi_{i-1}\in \B^G$ and
    $2^{i-1}\alpha_{i-1}=\alpha$.  By construction
    $(\beta_{i-1})_\p = 0$ at all archimedean $\p$, and so
    $\beta_{i-1}\in D(\Br(l))$. Now $\h_{\X^G}(\chi_{i-1})\ge
    \omega$, and hence we may apply Lemma~\ref{le:bggei} again to
    show that $\h_{\B^G}(\alpha_{i-1}) \ge \omega$.  Hence
    $\h_{\B^G}(\alpha)\ge \omega + (i-1)$.  If $h=\omega+n$ for
    some $n$, we have that $\omega+(n-1)=\h_{\X^G}(\chi)-1\le
    \h_{\B^G}(\alpha)$.  If $h=\omega 2$, we have shown that
    $\h_{\B^G}(\alpha)=\omega 2$. A similar argument handles the
    case $h=\infty$, where $\h_{\B^G}(\alpha)=\infty$ as well.
\end{proof}

\begin{lemma}\label{le:constructbeta}
    Suppose that $k$ is a number field, $\chi\in \X^G$, and
    $\h_{\X^G}(\chi)>1$. Then there exists $\alpha=\beta+\chi\in
    \B^G$ such that $\vert \alpha \vert = \vert \chi \vert$ and
    $\h_{\X^G}(\chi)-1\le \h_{\B^G}(\alpha) \le \h_{\X^G}(\chi)$.
\end{lemma}

\begin{proof}
    Since $\h_{\X^G}(\chi)>1$, $\chi=4\epsilon$ for some
    $\epsilon\in \X^G$. Let $\delta=2\epsilon$ and
    $\gamma'=s_1^*(\delta)$. Since $\h_{\X^G}(\delta)>0$,
    $\gamma'_{\p}=0$ at each archimedean place $\p$ of $l$.
    Similarly, since $\h_{\X^G}(\delta)>0$ and since by
    \eqref{eq:corsum2} the image $s_1^*(\X^G)$ is
    $\sigma$-negated, we have that $\gamma'_{\p}=0$ at every inert
    or ramified place $\p$ of $l$, and, by \eqref{eq:action}, that
    at every pair $\p_1, \p_2$ of $\sigma$-permuted places of $l$,
    $\gamma'_{\p_1} = -\gamma'_{\p_2}$ as well.

    We will define a $\gamma\in \Br(l)_2$ via its invariants.
    Let $\gamma_{\p}=0$ at every archimedean place of $l$.  For
    every pair $\p_1, \p_2$ of $\sigma$-permuted places of $l$,
    let $\gamma_{\p_1}=\gamma'_{\p_1}$ and $\gamma_{\p_2}=0.$  Now
    there are only finitely many such pairs $\p_1, \p_2$ at which
    $\gamma'$ has nontrivial invariants.  For each such pair
    $\p_1, \p_2$ choose an inert or ramified place $\q$ of $l$ and
    set $\gamma_{\q}=-\gamma_{\p_1}$.  At all other inert or
    ramified places $\q$ set $\gamma_\q=0$.

    By construction $\sum_\p \gamma_\p = 0$, so there exists a
    $\gamma\in \Br(l)_2$ with invariants $\{\gamma_\p\}$.
    Furthermore, $((1-\sigma)\gamma)_\p=0$ at all inert, ramified,
    or archimedean places $\p$ by \eqref{eq:action}. At pairs
    $\p_1, \p_2$ of $\sigma$-permuted places, we have
    \begin{equation*}
        ((1-\sigma)\gamma)_{\p_1}=\gamma'_{\p_1}=-\gamma'_{\p_2}=
        ((1-\sigma)\gamma)_{\p_2}.
    \end{equation*}
    Hence $(1-\sigma)\gamma = 2(1-\sigma)\gamma' = 2s_1^*(\delta) =
    s_1^*(2\delta) = s_1^*\chi$ and $\gamma+\delta\in
    \B^G$. Since the invariants of $\gamma$ are either 0 or equal
    to a corresponding invariant of $\gamma'$, which is a
    homomorphic image of $\delta$, $\vert \gamma \vert \le \vert
    \delta \vert$, and hence $\vert \gamma+\delta \vert = \vert
    \delta\vert$.  Setting $\beta=2\gamma$ and $\alpha=
    2(\gamma+\delta)= \beta+\chi$, we then have $\vert \alpha
    \vert = \vert \chi\vert$.  Moreover, $\h_{\B^G}(\alpha)>0$
    since $\alpha=2(\gamma+\delta)$.  Using
    Lemma~\ref{le:doublele}, we have that $\h_{\X^G}(\chi)-1\le
    \h_{\B^G}(\alpha) \le \h_{\X^G}(\chi)$.
\end{proof}

\begin{proof}[Proof of Proposition~\ref{pr:equiv1}]
    ($\Rightarrow$)  Let $\alpha=\beta + \chi$ be an order 2
    element in $\B^G$ with $\h_{\B^G}(\alpha)=\omega 2$.  By
    restriction, we have that $\h_{\B^G}(\alpha)\le
    \h_{\X^G}(\chi)$. Since the Ulm length of $\X$ is at most
    $\omega 2$ \cite[Thm.~1]{FS}, we have that $\h_{\X^G}(\chi)\in
    \{\omega 2, \infty\}$.

    Suppose that $\chi\in D(\X^G)$. By
    Lemma~\ref{le:constructbeta}, there exists $\beta'\in
    \Br(l)_2$ such that $\alpha'=\beta'+\chi\in \B^G$, $\vert
    \alpha'\vert = \vert \chi\vert$, and $\alpha'\in D(\B^G)$. But
    then $\gamma=\alpha-\alpha' = \beta-\beta'\in \B^G$ satisfies
    $\h_{\B^G}(\gamma)=\omega 2$. Now an element $b\in \Br(l)_2$
    lies in $\B^G$ if and only if $(1-\sigma)b=0$. But
    $\Br(l)_2^{\langle \sigma\rangle}$ consists, by
    \eqref{eq:action}, of a restricted direct sum of $\Z/2\Z$ and
    $\Q_2/\Z_2$ summands. Therefore there is no element in
    $\Br(l)_2\cap \B^G$ of height $\omega 2$ and we have a
    contradiction. Therefore $\h_{\X^G}(\chi)=\omega 2$.

    ($\Leftarrow$)  Now suppose that $\chi\in \X^G$ is an order 2
    element of height $\omega 2$.  By
    Lemma~\ref{le:constructbeta}, there exists a $\beta\in \Br(l)$
    such that $\alpha=\beta+\chi\in \B^G$ is of order 2 and
    $\h_{\B^G}(\alpha)=\omega 2$, and we are done.
\end{proof}

\subsection{Reducing to $\lambda$}\

\begin{proposition}\label{pr:order2iflambda}
    Suppose $k$ is a number field.  Then there exists an order 2
    element in $\B^G$ of height $\omega 2$ if and only if
    $\lambda_w=\omega 2$ for some nontrivial $w\in W_I$.
\end{proposition}

\begin{lemma}\label{le:kpinside}
    Suppose $k$ is a number field and $\sum\chi_p\in \Xn^G$ is an
    order 2 element with height greater than the Ulm length of
    $\X(l)_2$.  Then $\Cor \sum \chi_p\in W$.
\end{lemma}

\begin{proof}
    If $\Cor\sum\chi_p=0$, then $0\in W$ and we are done.

    Otherwise, let $w=\Cor\sum\chi_p$.  The order of $w$ is 2,
    and, since by \eqref{eq:corsum} $\Cor$ is a homomorphism
    from $\Xn^G$ to the $\sigma$-invariant subgroup of $\X(l)$,
    $w\in D(\X(l)_2^{\langle \sigma\rangle})\subset
    D(\X(l)_2)^{\langle \sigma\rangle}$.  Hence $w\in I\oplus
    {}_2N\oplus P^{\langle\sigma\rangle}$.

    Consider $p\neq \tilde p$ for which $\Cor \chi_p$ is not trivial.
    By \eqref{eq:corsigma},
    \begin{equation*}
        \Cor (\chi_p + \chi_{\tilde p}) = \Cor (\chi_p +
        s_{p\tilde p}^*(\chi_p)) = (1+\sigma)\Cor \chi_p.
    \end{equation*}
    Now $\Cor\chi_p$ is an element of order at most two in
    $\X(l)_2$, therefore represented by $l(\sqrt{e})$ for $e\in
    l^*$.  Then $(1+\sigma)\Cor \chi_p$ is represented by
    $l(\sqrt{N_{l/k}(e)})$.  Set $z_{p\tilde p}=N_{l/k}(e)$.

    Now consider $p=\tilde p$ for which $\Cor \chi_p$ is not
    trivial.  Then as in the proof of \eqref{eq:spppistp}, there
    exists a $\tau\in G_l$ such that $a_p=c/\sigma\tau(a_p)$ and
    we let $\tilde \sigma = \sigma \tau$. Then $\tilde \sigma$ is
    an automorphism of $l(a_p)$ of order 2 with fixed field $k_p
    := k(a_p+c/a_p)$.

    We claim that $[k_p:k]$ is even.  Since $\sqrt{d}\notin k_p$
    and $[l(a_p):k_p]=2$, $l(a_p)= k_p(\sqrt{d})$.  Then $a_p\in
    l(a_p)$ satisfies $N_{l(a_p)/k_p} a_p=c$, so the quaternion
    algebra $(c,d)$ splits over $k_p$.  But then $[k_p\colon k]$
    is even.

    Since $\chi_p\in \X(l(a_p))$ is of order 2, it is represented
    by $l(a_p)(\sqrt{f})$, where $f$ is determined up to its class
    in $l(a_p)^*/{l(a_p)^{*2}}$.  Moreover, since
    $s_{pp}^*(\chi_p)= \chi_p$, by \eqref{eq:spppistp},
    $\tilde\sigma(\chi_p)= \chi_p$, or $\tilde\sigma(f)=f$ in
    $l(a_p)^*/{l(a_p)^{*2}}$.  Hence $N_{l(a_p)/k_p}(f) \in
    l(a_p)^{*2}$.  By Kummer theory, $(k_p^*\cap
    l(a_p)^{*2})/k_p^{*2}$ consists only of the classes $1$
    and $d$.  Therefore, modulo $k_p^{*2}$, $N_{l(a_p)/k_p}(f)$ is
    either $1$ or $d$.

    Now $\Cor \chi_p$ is represented by $l(\sqrt{e})$, where
    $e=N_{l(a_p)/l}(f)$, and then $N_{l/k}(e)=N_{k_p/k}
    N_{l(a_p)/k_p}(f)$.  Since $N_{l(a_p)/k_p}(f)$ is $1$ or $d$
    mod $k_p^{*2}$ and $N_{k_p/k}(dz^2) = d^{[k_p:k]}
    N_{k_p/k}(z)^2 \in k^{*2}$, we deduce that $N_{l/k}(e)$ is a
    square in $k^*$. Then, from the square-class exact sequence
    (\cite[Thm.~3.4]{La})
    \begin{equation*}
        1\to \langle d\cdot {k^*}^2\rangle \to k^*/{k^*}^2 \to
        l^*/{l^*}^2 \xrightarrow{N_{l/k}} k^*/{k^*}^2
    \end{equation*}
    we have that, up to squares in $l^*$, $e$ is represented by an
    element of $k^*$. Set $z_{pp}$ to be this value.

    For all remaining $p$, set $z_{p\tilde p}=1$.

    Now $\Cor\sum\chi_p$ is represented by $l(\sqrt{e})$, where
    $e$ is the product in $l^*$ of $z_{p\tilde p}$ for all
    $\{p,\tilde p\}$. Hence $\Cor\sum\chi_p$ is an element of
    order at most 2 represented by $l(\sqrt{e})$, where $e$ is a
    product of elements from $k$.
\end{proof}

\begin{lemma}\label{le:lambdaomega2}
    Suppose $k$ is a number field and $\chi\in \X^G$ is an order 2
    element with $\h_{\X^G}(\chi)= \omega 2$. Then there exists
    $\hat\chi \in D(\X^G)$ with $\vert\hat\chi\vert=2$ such that
    $w=\chi-\hat\chi\in W_I$ is nontrivial and $\lambda_w=\omega
    2$.
\end{lemma}

\begin{proof}
    Write $\chi=\chi_u + \chi_{\neq u}$ with $\chi_{\neq u} = \sum
    \chi_p$. Since for each $p$ the Ulm length of $\X(l(a_p))_2$
    is $\omega+n$ for some $n\in \N_{\ge 0}$ \cite[Thm.~1]{FS} and
    since $\h_{\X^G}(\chi)=\omega 2$, we have $\chi_{\neq u}\in
    D(\Xn)$. For each summand $\X(l(a_p))_2$ with $p=\tilde p$,
    and for each pair of summands $\X(l(a_p))_2\oplus
    \X(l(a_{\tilde p}))_2$ for $p\neq \tilde p$, the divisible
    subgroup is a finite direct sum of $\Q_2/\Z_2$ components
    \cite[Thm.~11.1.2]{NSW}.

    The fixed subgroup of an order 2 action on a finite direct sum
    of $\Q_2/\Z_2$ components is a direct sum of $\Z/2\Z$ and
    $\Q_2/\Z_2$ components.  Hence for $\tilde p=p$ the
    $s_{pp}^*$-fixed subgroup of $D(\X(l(a_p))_2)$, and when $\tilde
    p\neq p$, the $(s_{p\tilde p}+s_{\tilde pp})^*$-fixed subgroup
    of $D(\X(l(a_p))_2\oplus \X(l(a_{\tilde p}))_2)$, is a direct sum
    of $\Z/2\Z$ and $\Q_2/\Z_2$ components.  Since $\chi_{\neq
    u}\in 2\X^G_{\neq u}$, in each summand or pair of summands,
    then, the components of $\chi_{\neq u}$ lie in the divisible
    part of the $\sum_p s_{p\tilde p}^*$-fixed subgroup. Hence
    $\chi_{\neq u}\in D(\Xn^G)$. For each $p$, let
    $\{\chi_p^{(i)}\}_{i=0}^\infty$ be a divisible tower over
    $\chi_p$, so that $\{\sum\chi_p^{(i)}\}_{i=0}^\infty \subset
    X^G_{\neq u}$ is a divisible tower over $\chi_{\neq u}$.

    Since $\h_{\X^G}(\chi)=\omega 2$ and the Ulm length of
    $\X(l)_2$ is $\omega+n$ for some $n\in \N_{\ge 0}$, $\chi_u\in
    D(\X(l)_2)$. Following the decomposition of $D(\X(l)_2)$ in
    \eqref{eq:dchil2}, write $\chi_u=w_I+w_N+w_P$.

    For each pair $\Q_2/\Z_2\oplus \Q_2/\Z_2$ of $\sigma$-permuted
    summands in $P$, denote the components of $z\in P$ in these
    summands by $z_s$ and $z_t$.  Define
    \begin{equation*}
        (w_P^{(i)})_s = - (\Cor \sum \chi_p^{(i)})_s -
        \frac{1}{2^i} (w_P)_t
    \end{equation*}
    and $(w_P^{(i)})_t = \frac{1}{2^i} (w_P)_t$. (We denote by
    $(1/2^i)(w_P)_t$ some element yielding $(w_P)_t$ under
    multiplication by $2^i$, and we fix this element for the
    duration.) Then $\vert w_P^{(i)}\vert \le 2^{i+1}$ and
    $\{w_P^{(i)}\}_{i=0}^\infty$ is a divisible tower over
    $w_P^{(0)}=w_P$ since $(1+\sigma)w_P = (-\Cor \sum \chi_p)_P$
    implies $((1+\sigma)w_P)_s = (-\Cor \sum \chi_p)_s$ and
    $(w_P)_s = (-\Cor \sum \chi_p)_s - (w_P)_t$.

    Now
    \begin{equation*}
        ((1+\sigma)w_P^{(i)})_s = - (\Cor \sum
        \chi_p^{(i)})_s = ((1+\sigma)w_P^{(i)})_t.
    \end{equation*}
    Furthermore, by \eqref{eq:corsum} the image of $\Cor$ on the
    divisible tower $\{\sum\chi_p^{(i)}\}_{i=0}^\infty$ over
    $\chi_{\neq u}=\sum \chi_p$ lies in the $\sigma$-invariant
    part of $P$, hence with components lying in the diagonals of
    the pairs $\Q_2/\Z_2\oplus \Q_2/\Z_2$.  Hence
    $((1+\sigma)w_P^{(i)})_t = -(\Cor \sum \chi_p^{(i)})_t$ as
    well.

    Since $N$ is divisible, we may choose a divisible tower
    $\{w_N^{(i)}\}_{i=0}^\infty\subset N$ over $w_N$. Now
    $(1+\sigma)w_N^{(i)}=0$ since $\sigma$ acts by negation on
    $N$. Furthermore, by \eqref{eq:corsum} the image of $\Cor$ on
    the divisible tower over $\chi_{\neq u}$ lies in $D(\X(l)_2)$
    and is $\sigma$-invariant, hence has zero component in $N$.
    Hence $(1+\sigma)(w_N^{(i)}+w_P^{(i)})= - (\Cor \sum
    \chi_p^{(i)})_{N\oplus P}$.

    Finally set $\hat w_I^{(i)} = (-\Cor\sum \chi_p^{(i+1)})_I$.
    Then since $\sigma$ is invariant on $I$,
    \begin{equation*}
        (1+\sigma)(\hat w_I^{(i)})= 2\hat w_I^{(i)}=2(-\Cor\sum
        \chi_p^{(i+1)})_I=(-\Cor \sum\chi_p^{(i)})_I.
    \end{equation*}
    Moreover, $\{\hat w_I^{(i)}\}_{i=0}^\infty$ is a divisible
    tower over $\hat w_I := \hat w_I^{(0)}$.  Note that $\hat w_I$
    has order at most 2 because $(-\Cor \sum \chi_p)_I =
    ((1+\sigma)\chi_u)_I = (1+\sigma)w_I = 0$ since $\chi\in \X^G$
    and $w_I$ is of order 2.

    Let $\hat\chi_u^{(i)}=\hat w_I^{(i)}+w_N^{(i)}+w_P^{(i)}$ and
    $\hat\chi^{(i)}=\hat\chi_u^{(i)}+\sum \chi_p^{(i)}$. Then
    $\{\hat\chi^{(i)}\}_{i=0}^\infty \subset \X^G$ is a divisible
    tower over $\hat\chi := \hat\chi_u^{(0)} + \sum \chi_p$.
    Clearly $w := \chi-\hat \chi=w_I-\hat w_I\in I$ and $\vert
    w\vert \le 2$. If $w=0$ then $\chi=\hat\chi$ and we have a
    contradiction: $\chi\in D(\X^G)$.  Hence $\vert w\vert = 2$
    and $\h_{\X^G}(w)=\omega 2$.

    Now since $\h_{\X^G}(w)=\omega 2$, for any $n\in \N$ there
    exists a $\chi'\in \X^G$ of height $\omega+n$ and $w=2\chi'$.
    We restrict $\omega+n$ to ordinals greater than the Ulm length
    of $\X(l)_2$. Write $\chi'=\chi'_u+\sum\chi'_p$.  Then
    $\chi'_u\in D(\X(l)_2)$.  Moreover, $(1+\sigma)\chi'_u=-\Cor
    \sum \chi'_p$. Now $\vert \sum\chi'_p \vert$ is at most 2
    since $w=2\chi'$ has no component in $\Xn$. If the order is 2,
    then by Lemma~\ref{le:kpinside}, $-\Cor\sum\chi'_p$ lies
    in $W$; if the order is 1, then $-\Cor\sum\chi'_p=0\in W$.

    Since $\chi'_u\in D(\X(l)_2)$, write $\chi'_u =
    \chi'_I+\chi'_N+\chi'_P$ according to \eqref{eq:dchil2}.  Since
    $2\chi'=w\in I$, $\chi'_N$ and $\chi'_P$ are of order at most
    2.  Hence $(1+\sigma)\chi'_I= 2\chi'_I=w_I=w$ and
    $(1+\sigma)\chi'_N= 2\chi'_N=0$. Therefore
    $(1+\sigma)\chi'_u=-\Cor \sum \chi'_p = w+(1+\sigma)\chi'_P\in
    W$. Therefore $(-\Cor \sum \chi'_p)_I = w$, and
    $\h_{\Xn^G}(\sum\chi'_p)\ge \h_{\X^G}(\chi')=\omega+n$. Hence,
    by considering $-\sum\chi'_p$ instead of $\sum\chi'_p$, we
    have $\lambda_w\ge \omega 2$.

    Now suppose that $\lambda_w=\infty$.  Then since the Ulm
    length of $\X^G_{\neq u}$ is at most $\omega 2$
    \cite[Thm.~1]{FS}, there exists $\chi_{\neq u}=\sum \chi_p\in
    \Xn^G$ of order at most 2 with $\h_{\Xn^G}(\chi_{\neq
    u})=\infty$ and $(\Cor\chi_{\neq u})_I=w$. Let $\{\sum
    \chi_p^{(i)}\}_{i=0}^\infty \subset \X^G_{\neq u}$ be a
    divisible tower over $\chi_{\neq u}$.  Now let $\hat P$ be a
    finite direct sum of $\Q_2/\Z_2$ summands in $P$ such that
    $\hat P \oplus \sigma(\hat P) = P$. For each $i\ge 1$, set
    $\chi_u^{(i)}=(-\Cor\sum\chi_p^{(i+1)})_I +
    (-\Cor\sum\chi_p^{(i)})_{\hat P}$.  Then
    \begin{equation*}
        ((1+\sigma)(\chi_u^{(i)}))_I = 2(\chi_u^{(i)})_I=
        (-2\Cor\sum\chi_p^{(i+1)})_I= (-\Cor\sum\chi_p^{(i)})_I,
    \end{equation*}
    and
    \begin{equation*}
        ((1+\sigma)(\chi_u^{(i)}))_{\hat P} = (\chi_u^{(i)})_{\hat
        P}= (-\Cor\sum\chi_p^{(i)})_{\hat P}.
    \end{equation*}
    Moreover, since by \eqref{eq:corsum} the image of $\Cor$ is
    $\sigma$-invariant on $\Xn^G$, the equality holds over
    $\sigma(\hat P)$ as well, and holds over $N$ since both sides
    must be trivial on $N$. Then for $i\ge 0$, $\chi_u^{(i)}+
    \sum\chi_p^{(i)}\in \X^G$.  Now
    \begin{align*}
        2\chi_u^{(0)} &= 2((-\Cor\sum\chi_p^{(1)})_I+
        (-\Cor\sum\chi_p^{(0)})_{\hat P})\\
        &=(-\Cor\sum\chi_p^{(0)})_I=-w.
    \end{align*}
    Since $\vert w\vert \le 2$, $-w=w$. Let $\chi'_0=w$ and for
    $i\ge 1$, $\chi'_i=\chi_u^{(i-1)}+\sum\chi_p^{(i-1)}$.  Then
    $2\chi'_1=w$ and $\{\chi'_i\}_{i=0}^\infty \subset \X^G$ is a
    divisible tower over $w$, so $w\in D(\X^G)$. But then
    $\chi=w+\hat\chi\in D(\X^G)$, a contradiction.
\end{proof}

\begin{proof}[Proof of Proposition~\ref{pr:order2iflambda}]
    ($\Rightarrow$) By Proposition~\ref{pr:equiv1}, if there is an
    order 2 element in $\B^G$ of height $\omega 2$, then there is
    an order 2 element in $\X^G$ of height $\omega 2$.  By
    Lemma~\ref{le:lambdaomega2}, there exists a nontrivial $w\in W_I$
    with $\lambda_w=\omega 2$.

    ($\Leftarrow$) Suppose that $\lambda_w=\omega 2$ for some
    nontrivial $w\in W_I$. We claim that there is an order 2
    element $\hat w$ in $\X^G$ with $\h_{\X^G}(\hat w)=\omega 2$.

    Since $\lambda_w=\omega 2$, for each ordinal $\omega+i$
    greater than the Ulm length of $\X(l)_2$ there exists an
    element $\chi_{\neq u}=\sum\chi_p\in \Xn^G$ with
    $\h_{\Xn^G}(\chi_{\neq u})=\omega+i$ and $-(\Cor \chi_{\neq
    u})_I=-w=w$.

    For each such ordinal $\omega+i$, we proceed as follows.  Let
    $\chi_{\neq u}^{(i)}\in \Xn^G$ be an element such that
    $2^i\chi_{\neq u}^{(i)}=\chi_{\neq u}$ and for each $j\in \N$
    let $\delta_{\neq u}^{(j)}\in \Xn^G$ satisfy $2^j\delta_{\neq
    u}^{(j)}=\chi_{\neq u}^{(i)}$.  Now set $\chi_u^{(i)}=-\Cor
    \chi_{\neq u}^{(i)}$ and $\chi^{(i)}=\chi_u^{(i)}+\chi_{\neq
    u}^{(i-1)}$.  Then since by \eqref{eq:corsum} the image of
    $-\Cor$ on $\Xn^G$ is $\sigma$-invariant,
    \begin{equation*}
        (1+\sigma)\chi_u^{(i)}=2\chi_u^{(i)}=-2 \Cor \chi_{\neq
        u}^{(i)}=-\Cor \chi_{\neq u}^{(i-1)}
    \end{equation*}
    and we have $\chi^{(i)}\in \X^G$.  Moreover, $2^i
    \chi^{(i)} = -\Cor \chi_{\neq u}$ with $-(\Cor \chi_{\neq
    u})_I = w$. Continuing on, for each $j\ge 2$, set
    $\delta_u^{(j)}=-\Cor \delta_{\neq u}^{(j)}$ and
    $\delta^{(j)}=\delta_u^{(j)}+ 2\delta_{\neq u}^{(j)}$; for
    $j=1$ set $\delta_u^{(1)}=-\Cor \delta_{\neq u}^{(1)}$ and
    $\delta^{(1)}=\delta_u^{(1)}+ \chi_{\neq u}^{(i)}$.  Then
    $2^j\delta^{(j)}=\chi^{(i)}$ and as before, $\delta^{(j)}\in
    \X^G$ for each $j$.  Hence we have shown that $w'=-\Cor
    \chi_{\neq u}$ is an order 2 element in
    $\X(l)^{\langle\sigma\rangle}$ with $\h_{\X^G}(w')\ge
    \omega+i$ and $w'_I=w$. In fact, because the height is greater
    than the Ulm length of $\X(l)_2$, $w'\in
    D(\X(l)_2)^{\langle\sigma\rangle}$.

    Hence for every $\omega+i$ we have produced an order 2 element
    $w'_i\in D(\X(l)_2)^{\langle\sigma\rangle}$ which satisfies
    $\h_{\X^G}(w'_i)\ge \omega+i$ and $(w'_i)_I=w$.  Now the
    divisible subgroup of $\X(l)_2$ is a finite direct sum of
    $\Q_2/\Z_2$ summands \cite[Thm.~11.1.2]{NSW}, and hence its
    exponent 2 subgroup is finite.  Hence for some $\hat w\in
    D(\X(l)_2)^{\langle \sigma\rangle}$, $\h_{\X^G}(\hat w)\ge
    \omega + i_n$ for an unbounded strictly increasing
    $\{i_n\}_{n=1}^\infty$ of natural numbers.  Therefore there
    exists a $\hat w$ of order 2 in
    $D(\X(l))^{\langle\sigma\rangle} \subset \X^G$ with
    $\h_{\X^G}(\hat w)\ge \omega 2$ and $\hat w_I=w$.

    Now suppose that $\hat w$ is divisible in $\X^G$.  Let
    $\{\chi^{(i)}\}_{i=0}^\infty \subset \X^G$ be a divisible
    tower over $\hat w$.  Write
    $\chi^{(i)}=\chi_u^{(i)}+\chi_{\neq u}^{(i)}$ for each $i$,
    and consider $\chi'=\chi_{\neq u}^{(1)}$, necessarily of order
    less than or equal to 2 since $\chi_{\neq u}^{(0)}=0$ because
    $\hat w\in \X(l)_2$. Then $\{\chi_{\neq
    u}^{(i+1)}\}_{i=0}^\infty \subset \Xn^G$ is a divisible tower
    over $\chi'$.  Since $\chi^{(1)}\in \X^G$, we have that
    $(1+\sigma)\chi_u^{(1)}=-\Cor\chi_{\neq u}^{(1)}$. Restricting
    the equation to $I$, we have $(\chi_u^{(0)})_I =
    (-\Cor\chi')_I$, and the left hand side is in fact $(\hat
    w)_I=w$.  Hence $\chi'$ is an order 2 element in $\Xn^G$ with
    $(\Cor\chi')_I=w$; we have then that $\lambda_{w}=\infty$, a
    contradiction.
\end{proof}

\section{Proof of Main Theorem}\label{se:proof}

\begin{proof}
    We show that any order two element $\alpha$ in $\Br(E)$ of
    height at least $\omega 2$ is divisible.

    Suppose $\alpha$ is an order two element in $\Br(E) =
    H^{2}(G_{k},\bar k(u)^{*})$ of height at least $\omega 2$.
    Then the image $\phi(\alpha)$ of $\alpha$ in
    \eqref{eq:brexactseq} is an element of $\B^G$ of height at
    least $\omega 2$.

    Since the Leopoldt conjecture holds for $l$ and $2$, there is
    only one $\Z_2$-extension of $l$.  Since $l$ is the quadratic
    subextension of $k^{cyc}/k$, the $\Z_2$-extension of $l$ is
    precisely $k^{cyc}$.  Hence $D(\X(l)_2)\cong \Q_2/\Z_2$ and
    $\sigma$ acts trivially on $D(\X(l)_2)$.  In the decomposition
    of $D(\X(l)_2)$ in \eqref{eq:dchil2}, $N=P=0$ and $I\cong
    \Q_2/\Z_2$. But since the quadratic subextension of
    $k^{cyc}/l$ is cyclic of order 4 over $k$, $W=0$.  Hence by
    Proposition~\ref{pr:order2iflambda}, there exists no order 2
    element of $\B^G$ of height $\omega 2$.  Hence $\phi(\alpha)$
    is divisible.

    If $\phi(\alpha)\neq 0$ and $\{\tilde \alpha_i\}$ is a
    divisible tower over $\phi(\alpha)$, then we may find
    preimages $\alpha_{n}\in \Br(E)$ such that
    $\phi(\alpha_{n})=\tilde\alpha_{n}$ and $2^{m}\alpha_{n} =
    \alpha_{n-m}$, as follows.  Let $\alpha_{n}=2\cdot
    \phi^{-1}(\tilde\alpha_{n+1})$.  Since the kernel of $\phi$ is
    the relative Brauer group $\Br(El/E)$ and is therefore of
    exponent 2, this map is well-defined.  We calculate
    \begin{equation*}
        \phi(\alpha_{n})=\phi(2\cdot \phi^{-1}(\tilde
        \alpha_{n+1}))= 2 \cdot \phi(\phi^{-1}(\tilde
        \alpha_{n+1}))=2 \cdot \tilde \alpha_{n+1}=\tilde
        \alpha_{n}
    \end{equation*}
    and
    \begin{align*}
        2^{m}\alpha_{n} & =2^{m}\cdot 2\cdot
        \phi^{-1}(\tilde\alpha_{n+1}) =2\cdot 2^{m}\cdot
        \phi^{-1}(\tilde\alpha_{n+1}) \\ &=2\cdot
        \phi^{-1}(2^{m}\tilde\alpha_{n+1})= 2\cdot
        \phi^{-1}(\tilde\alpha_{n-m+1}) =\alpha_{n-m}.
    \end{align*}
    Hence $\{\alpha_n\}$ is a divisible tower over $\alpha$ and
    $\alpha$ is a divisible element of $\Br(E)$.

    If $\phi(\alpha)=0$, $\alpha$ lies in the kernel of $\phi$ and
    so is split by base change to $El=l(u)$.  Hence $\alpha$ is
    represented by the class of a central simple algebra with
    $E(\sqrt{d})/E$ as a maximal subfield.  Suppose that
    $\alpha=\epsilon_1$ for a quaternion algebra $\epsilon_{1}=
    (d,e)_{E}$ with $e\in E^{\times}$.  Now since $k^{cyc}$ and
    $E$ are linearly disjoint over $k$, $k^{cyc}E$ is a
    $\Z_{2}$-extension of $E$ containing $El$.  Let $k^{cyc}_nE$
    be the $\Z/2^{n}\Z$ layer of $k^{cyc}E/E$, and choose
    generators $\sigma_{n}\in G_{k^{cyc}_{n}E/E}$ satisfying
    ${\sigma_{n}}\vert_{k_{n-1}^{cyc}E} = \sigma_{n-1}$.  Then the
    cyclic algebras $\epsilon_{n} = (k^{cyc}_{n}E/E,\langle
    \sigma_{n}\rangle, e)$ are each of order $2^{n}$ in $\Br(E)$
    and, moreover, $2^{m}\epsilon_{n}=\epsilon_{n-m}$.  Hence
    $\alpha$ is divisible in $\Br(E)$.

    We have shown that $U_2(\omega 2, \Br(E)) = 0$.  By
    \cite{FSS1} and \cite{FSS2}, $U_2(\omega 2, \Br(k(t))) = 0$
    and all other Ulm invariants of $\Br(E)$ and $\Br(k(t))$ are
    identical. Therefore $\Br(E)\cong \Br(k(t))$.
\end{proof}

\section*{Acknowledgments}

The second author thanks the Department of Mathematics at the
Technion---Israel Institute of Technology for its hospitality
during 1998--1999.


\begin{thebibliography}{NSW}
\bibliographystyle{alpha}

\bibitem[AB]{AB}
M.~Auslander and A.~Brumer.  Brauer groups of discrete valuation
rings.  \textit{Nederl.~Akad.~Wetensch.~Proc.~Ser.~A} \textbf{30}
(1968), 286--296.

\bibitem[Br]{Br}
A.~Brumer.  On the units of algebraic number fields.
\textit{Mathematika} \textbf{14} (1967), 121--124.

\bibitem[CR]{CR}
C.~Curtis and I.~Reiner.  \textit{Representation theory of finite
groups and associative algebras}.  Pure and Applied Mathematics
11.  New York: Wiley-Interscience, 1962.

\bibitem[Fa]{F}
D.~K.~Faddeev.  Simple algebras over a field of algebraic
functions of one variable. \textit{Trudy Mat.~Inst.~Steklov}
\textbf{38} (1951), 321--344; \textit{Amer.~Math.~Soc.~Transl.} II
\textbf{3} (1956), 15--38.

\bibitem[FS]{FS}
B.~Fein and M.~Schacher.  Brauer groups and character groups of
function fields.  \textit{J.~Algebra} \textbf{61} (1979),
249--255.

\bibitem[FSS1]{FSS1} B.~Fein, M.~Schacher, and J.~Sonn.
Brauer groups of rational function fields.  \textit{Bull. Amer.
Math. Soc. (N.S.)} \textbf{1} (1979), no. 5, 766--768.

\bibitem[FSS2]{FSS2}
B.~Fein, M.~Schacher, and J.~Sonn.  Brauer groups of fields of
genus zero.  \textit{J.~Algebra} \textbf{114} (1988), no.~2,
479--483.

\bibitem[Ja]{Ja}
G.~Janusz.  Automorphism groups of simple algebras and group
algebras. \textit{Representation theory of algebras (Proc.~Conf.,
Temple Univ., Philadelphia, Pa., 1976)}. Lecture Notes in Pure
Appl.~Math. 37. New York: Dekker, 1978, pp.~381--388.

\bibitem[La]{La} T.~Y.~Lam.  \textit{The algebraic theory of
quadratic forms}, revised 2nd printing.  Mathematics Lecture Note
Series. Reading, Mass.: Benjamin/Cummings Publishing~Co., Inc.,
1980.

\bibitem[NSW]{NSW}
J.~Neukirch, A.~Schmidt, and K.~Wingberg.  \textit{Cohomology of
number fields}.  Grundlehren der mathematischen Wissenschaften
323. Berlin: Springer-Verlag, 2000.

\end{thebibliography}
\end{document}